\documentclass[12pt,reqno]{amsart}

\usepackage{amssymb,enumerate}

\newcommand{\ds}{\displaystyle}     

\DeclareMathAlphabet{\E}{U}{eus}{m}{n}     

\newcommand{\V}{{\mathcal V}}

\newcommand{\PP}{{\mathbb P}}

\newcommand{\N}{{\mathbb N}}

\newcommand{\kk}{{\Bbbk}}

\newcommand{\la}{\langle}
\newcommand{\ra}{\rangle}

\newtheorem{thm}{Theorem}[section]
\newtheorem{lemma}[thm]{Lemma}

\newtheorem{prop}[thm]{Proposition}

\theoremstyle{definition}
\newtheorem{defn}[thm]{Definition}
         
\newtheorem{notn}[thm]{Notation}

\newtheorem{rmk}[thm]{Remark}

\newtheorem{example}[thm]{Example}

\newtheorem{Pf}{Proof$\!\!$}         
         
\newenvironment{pf}{\begin{Pf}}{\qed\end{Pf}}

\DeclareMathSymbol{\twoheadrightarrow}  {\mathrel}{AMSa}{"10}

\newcounter{letter}
\renewcommand{\theletter}{\rom{(}\alph{letter}\rom{)}}

\newcounter{rnum}
\renewcommand{\thernum}{\rom{(}\roman{rnum}\rom{)}}

\begin{document}
\baselineskip21pt


\title[Point Modules over Regular Graded Skew Clifford Algebras]%
{Point Modules over \\[3mm] Regular Graded Skew Clifford Algebras}


\maketitle

\vspace*{0.1in}

\baselineskip15pt

\renewcommand{\thefootnote}{\fnsymbol{footnote}}
\centerline{\sc Michaela Vancliff\footnote{The first author was
supported in part by NSF grants DMS-0457022 \& DMS-0900239.\\[-3mm]}}
\centerline{Department of Mathematics, P.O.~Box 19408}
\centerline{University of Texas at Arlington,
Arlington, TX 76019-0408}
\centerline{{\sf vancliff@uta.edu}}
\centerline{{\sf www.uta.edu/math/vancliff}}

\bigskip
\centerline{and}
\bigskip

\centerline{\sc Padmini P. Veerapen\footnote{%
\begin{minipage}[t]{6.6in}%
The second author was supported in part by NSF grant DGE-0841400 as a 
Graduate Teaching Fellow in U.T.\ Arlington's GK-12 MAVS Project.%
\end{minipage}%
}}
\centerline{Department of Mathematics, P.O.~Box 19408}
\centerline{University of Texas at Arlington,
Arlington, TX 76019-0408}
\centerline{{\sf pveerapen@uta.edu}}

\setcounter{page}{1}
\thispagestyle{empty}

\bigskip
\bigskip

\begin{abstract}
\baselineskip15pt
Results of Vancliff, Van Rompay and Willaert in 1998 (\cite{VVW}) prove 
that point modules over a regular graded Clifford algebra (GCA) are 
determined by (commutative) quadrics of rank at most two that
belong to the quadric system associated to the GCA.  In 2010, in
\cite{CV}, Cassidy and Vancliff generalized the notion of a GCA to that 
of a graded skew Clifford algebra (GSCA). The results in this article 
show that the results of \cite{VVW} may be extended, with suitable
modification, to GSCAs. In particular, using the notion of $\mu$-rank 
introduced recently by the authors in \cite{VV}, the point modules over 
a regular
GSCA are determined by (noncommutative) quadrics of $\mu$-rank at most 
two that belong to the noncommutative quadric system associated to the 
GSCA.  
\end{abstract}

\baselineskip18pt

\newpage


\section*{Introduction}

The notion of a graded skew Clifford algebra (GSCA) was introduced in 
\cite{CV}, and it is an algebra that may be viewed as a quantized analog 
of a graded Clifford algebra (GCA). In \cite{CV}, it was shown that
many of the results that hold for GCAs have analogous counterparts in the
context of GSCAs.  In particular, homological and algebraic properties of
a GSCA, $A$, are determined by properties of a certain quadric system
associated to $A$.  The importance of GSCAs was
highlighted in \cite{NVZ}, where they were shown to play a critical
role in the classification of the quadratic AS-regular algebras of global
dimension three. Hence, GSCAs are expected to play a critical role in
the classification of the quadratic AS-regular algebras of global 
dimension four and greater. The reader is referred to \cite{AS,ATV1,ATV2}
for results concerning AS-regular algebras and their associated geometric
data, and to \cite{SV,VVW} for results concerning GCAs and their 
associated geometric data.

Consequently, it is reasonable to attempt to extend the results 
in~\cite{VVW} concerning point modules over GCAs to point modules over 
GSCAs. As such, our main objective in this article is to generalize 
\cite[Theorem~1.7]{VVW}. That result states, in part, that if the number,
$N$, of point modules over a regular GCA, $C$, is finite, then 
$N = 2 r_2 + r_1$,
where $r_j$ is the number of elements of rank $j$ that belong to the
projectivization of a certain quadric system associated to $C$ (see
Theorem~\ref{vvw} for the precise statement). We
achieve our objective in Theorem~\ref{bigthm}, where the notion of
$\mu$-rank (introduced in \cite{VV}) is used in place of the traditional 
notion of rank. However, 
\cite[Theorem~1.7]{VVW} also states, in part, that if $N < \infty$,
then $r_1 \in \{0,\ 1\}$. We present examples in the last section that 
demonstrate that this part of \cite[Theorem~1.7]{VVW} appears not to have
an obvious counterpart in the setting of GSCAs.

Although the flow of this article follows that of \cite[\S1]{VVW}, many
of our results require methods of proof that differ substantially from 
those used in \cite[\S1]{VVW}, since the proofs in \cite{VVW} make use 
of standard results concerning symmetric matrices and the general linear 
group. This article consists of two sections: in Section~\ref{sec1}, 
notation and terminology are defined, while 
Section~\ref{sec2} is devoted to proving our main result, which is
given in Theorem~\ref{bigthm}. 


\bigskip

\section{Graded Skew Clifford Algebras } \label{sec1}

In this section, we define the notion of a graded skew Clifford algebra 
from \cite{CV}, and give the relevant results from \cite{CV}
needed in Section~\ref{sec2}.

Throughout the article, $\kk$ denotes an algebraically closed field such 
that char$(\kk)\neq~2$, and $M(n,\ \kk)$ denotes the vector space of 
$n \times n$ matrices with entries in $\kk$. For a graded 
$\kk$-algebra~$B$, the span of the homogeneous elements in $B$ of 
degree~$i$ will be denoted $B_i$, and $T(V)$ will denote the
tensor algebra on the vector space~$V$.  If $C$ is a vector space, then 
$C^\times$ will denote the nonzero elements in $C$, and $C^*$ will denote 
the vector space dual of $C$. 

For $\{i, j\} \subset \{1, \ldots , n\}$, let $\mu_{ij} \in \kk^{\times}$ 
satisfy the property that $\mu_{ij}\mu_{ji} = 1$ for all $i,\ j$
where $i \ne j$. We write $\mu = (\mu_{ij}) \in M(n,\ \kk)$. As in 
\cite{CV}, we write $S$ for the quadratic $\kk$-algebra on generators 
$z_1, \ldots, z_n$ with defining relations $z_jz_i = \mu_{ij}z_iz_j$ 
for all $i, j = 1, \ldots , n$, where $\mu_{ii} = 1$ for all $i$. We set 
$U \subset T(S_1)_2$ to be the span of the defining relations of $S$
and write $z = (z_1, \ldots, z_n)^T$.

\begin{defn}\cite[\S1.2]{CV}\label{musymm}
\begin{enumerate}
\item[(a)]
With $\mu$ and $S$ as above, a (noncommutative) quadratic form is
defined to be any element of $S_2$.
\item[(b)]
A matrix $M \in M(n,\ \kk)$ is called $\mu$-symmetric if $M_{ij} = 
\mu_{ij}M_{ji}$ for all $i$, $j = 1, \ldots , n$. 
\end{enumerate}
\end{defn}

Henceforth, we assume $\mu_{ii} = 1$ for all $i$, and write 
$M^{\mu}(n,\ \kk)$ for the vector space of $\mu$-symmetric $n \times n$ 
matrices with entries in $\kk$. By \cite{CV}, there is a one-to-one
correspondence between elements of $M^{\mu}(n,\ \kk)$ and $S_2$ via 
$M \mapsto z^T M z \in S$.

\begin{notn}\label{tau}
Let $\tau : \PP(M^{\mu}(n,\ \kk)) \to \PP (S_2)$ be defined by 
$\tau (M) = z^TMz$.
\end{notn}

\begin{rmk}\label{linindep}
Henceforth, we fix $M_1, \ldots , M_n \in M^{\mu}(n,\ \kk)$. For
each $k = 1, \ldots , n$, we fix representatives $q_k = \tau(M_k)$. By 
\cite[Lemma 1.3]{CV}, $\{q_k\}_{k=1}^n$ is linearly independent in $S$ 
if and only if $\{M_k\}_{k=1}^n$ is linearly independent. This
correspondence mirrors the correspondence between symmetric matrices and 
commutative quadratic forms.
\end{rmk}

\begin{defn}\cite{CV}\label{GSCA}
A {\em graded skew Clifford algebra} $A = A(\mu, M_1, \ldots , M_n)$ 
associated to $\mu$ and $M_1$, $\ldots ,$ $M_n$ is a graded $\kk$-algebra
on degree-one generators $x_1, \ldots , x_n$ and on degree-two generators
$y_1, \ldots , y_n$ with defining relations given by:
\begin{enumerate}
\item[(a)] $\ds x_i x_j + \mu_{ij} x_j x_i = \sum_{k=1}^n (M_k)_{ij} y_k$
           for all $i, j = 1, \ldots , n$, and
\item[(b)] the existence of a normalizing sequence 
           $\{ y_1', \ldots , y_n'\}$ that spans 
	   $\kk y_1 + \cdots + \kk y_n$.
\end{enumerate}
\end{defn} 

\begin{rmk}\label{yirelns}
If A is a graded skew Clifford algebra (GSCA), then \cite[Lemma 1.13]{CV}
implies that $y_i \in (A_1)^2$ for all $i = 1, \ldots , n$ if and only if 
$M_1, \ldots , M_n$ are linearly independent. Thus, hereafter, we assume 
that $M_1, \ldots , M_n$ are linearly independent. 
\end{rmk}

By \cite{CV}, the degree of the defining relations of $A$ and certain
homological properties of $A$ are intimately tied to certain geometric 
data associated to $A$ as follows. 


\begin{defn}\cite{CV}\label{quadric}\hfill\\
\indent (a) \  
Let $\V(U) \subset \PP((S_1)^*) \times \PP((S_1)^*)$ denote the zero 
locus of $U$. For any $q \in S_2^{\times}$, we call the zero locus of 
$q$ in $\V(U)$ the {\em quadric associated to}~$q$, and denote it by 
$\V_{_U}(q)$; in other words, $\V_{_U}(q) = \V(\kk \hat q + U) = 
\V(\hat q) \cap \V(U)$, where $\hat q$ is any lift of $q$ to $T(S_1)_2$. 
The span of elements $Q_1, \ldots , Q_m\in S_2$ will be called the 
{\em quadric system} associated to $Q_1, \ldots , Q_m$.\\
\indent (b) \  
If a quadric system is given by a normalizing sequence in $S$, then it is
called a {\em normalizing quadric system}.\\
\indent (c) \  
We call a point $(a, \ b) \in \V(U)$ a {\em base point} of the quadric 
system associated to $Q_1, \ldots , Q_m \in S_2$ if  $(a, \ b) \in 
\V_{_U}(Q_k)$ for all $k = 1, \ldots , m$.
We say such a quadric system is {\em base-point free} if 
$\bigcap_{k=1}^m \V_{_U}(Q_k)$ is empty.
\end{defn}

\begin{thm}\label{main2}\cite{CV}
For all $k = 1, \ldots , n$, let $M_k$ and $q_k$ be as in 
Remark~\ref{linindep}. A graded skew Clifford algebra 
$A = A(\mu, M_1, \ldots , M_n)$ is a quadratic, Auslander-regular algebra
of global dimension $n$ that satisfies the Cohen-Macaulay property with 
Hilbert series $1/(1-t)^n$ if and only if the quadric system associated 
to $\{q_1 , \ldots , q_n\}$ is normalizing and base-point free; in this 
case, $A$ is a noetherian Artin-Schelter regular domain and is unique up 
to isomorphism.
\end{thm}

\begin{rmk}\label{W}\hfill\\
\indent (a) \  
Henceforth, we assume that the quadric system associated to 
$\{q_1 , \ldots , q_n\}$ is normalizing and base-point free. By 
Theorem~\ref{main2}, this assumption allows us to write 
$A = T(V)/\langle W \rangle$, where $V = (S_1)^*$ and $W \subseteq 
T(V)_2$,  and write the Koszul dual of $A$ as 
$T(S_1)/\langle W^{\perp} \rangle = S/\langle q_1, \ldots , q_n\rangle$. 
In this setting, $\{x_1, \ldots , x_n\}$ is the dual basis in $V$ to 
the basis $\{z_1, \ldots , z_n\}$ of $S_1$, and we write 
$\sum_{i, j} \alpha_{ijm} (x_ix_j + \mu_{ij}x_jx_i) = 0$ for the defining 
relations of $A$, where $\alpha_{ijm} \in \kk$ for all $i, j, m$, and 
$1 \le m \le n(n-1)/2$.\\
\indent (b) \  
By \cite[Lemma 5.1]{CV} and its proof, the set of pure tensors in
$\PP(W^{\perp})$, that is, $\{a \otimes b \in \PP(W^{\perp}): a, b
\in S_1 \}$, is in one-to-one correspondence with the zero locus
$\Gamma$, in 
$\PP(S_1) \times \PP(S_1)$, of W given by $\Gamma = \{(a,\ b) \in 
\PP(S_1) \times \PP(S_1) : w(a, b) = 0 \mbox{ for all } w \in W\}$. 
\end{rmk}

We will now make more precise the connection between points in the zero 
locus of $W$ and certain quadratic forms.

\begin{lemma}\label{qformzlocus}
If $a, b \in S_1^{\times}$, then the quadratic form $ab \in 
\PP\big(\sum_{i=1}^n \kk q_i\big)$ if and only if $(a,\ b) \in \Gamma$.
\end{lemma}
\begin{pf}
Suppose $w((a, b)) = 0$. By Remark~\ref{W}(b), $w(a \otimes b)
= 0$ for all $w \in W$, and so $a \otimes b \in W^{\perp}$. Since $S$ is 
a domain, $ab \ne 0$ in $S$, so $ab \in 
\PP\big(\sum_{i=1}^n \kk q_i\big)$, 
as desired. This argument is reversible, so the converse holds.
\end{pf}


\newpage

\section{Point Modules over Graded Skew Clifford Algebras}\label{sec2}

In this section, we prove results that relate point modules over 
GSCAs to noncommutative quadrics in the sense of Definition~\ref{quadric}.
In particular, we use the notion of $\mu$-rank introduced in \cite{VV} 
to extend results in \cite{VVW} about graded Clifford algebras (GCAs) to
GSCAs, with our main result being Theorem~\ref{bigthm}.
Although the overall approach and some of the proofs are 
influenced by those in \cite[\S1]{VVW}, many of the proofs involve new 
arguments. 

\bigskip

In \cite{VV}, a notion of $\mu$-rank of a (noncommutative) quadratic form
on $n$ generators was defined, where $n = 2$ or $3$. The results in 
\cite{VV} suggest a notion of $\mu$-rank at most two of a (noncommutative)
quadratic form on $n$ generators for any $n \in \N$ as follows.

\begin{defn}\cite{VV}\label{mu-rank}
Let $S$ be as in Definition~\ref{musymm}, where $n$ is an arbitrary
positive integer, and let $Q \in S_2$.
\begin{enumerate}
\item[{\rm (a)}]
If $Q = 0$, we define $\mu$-rank$(Q) = 0$.
\item[{\rm (b)}]
If $Q = L^2$ for some $L \in S_1^\times$, we define $\mu$-rank$(Q) = 1$.
\item[{\rm (c)}]
If $Q \neq L^2$ for any $L \in S_1^\times$, but $Q = L_1 L_2$ where
$L_1$, $L_2 \in S_1^\times$, we define $\mu$-rank$(Q) = 2$.
\end{enumerate}
Moreover, if $M \in \PP(M^{\mu}(n,\ \kk))$ and if $\mu$-rank$(\tau (M)) 
\leq 2$, where $\tau$ is given in Notation~\ref{tau}, then we define 
$\mu$-rank$(M)$ to be the $\mu$-rank of $\tau (M)$. 
\end{defn}

\begin{rmk}\label{quadTwoWays3}
In contrast to the commutative setting, there exist noncommutative 
quadratic forms $q$ where $0\neq q = L^2 = L_1 L_2$, with $L, L_1, L_2 
\in S_1$ and $L_1, \ L_2$ linearly independent. For example, let 
$n = 2 = \mu_{12}$ and $q = (z_1 + 2 z_2)^2 = (z_1 + z_2)(z_1 + 4 z_2)$.
\end{rmk}

We now define a function $\Phi$ that will play a role similar to that 
played by the function $\phi$ in \cite[Section~1]{VVW}.

\begin{defn}\label{Phi}
Let $a,\ b \in \PP^{n-1}$, with $a = (a_1, \ldots , a_n)$, 
$b = (b_1, \ldots , b_n)$, where $a_i$, $b_i \in \kk$ for all $i$.  
We define $\Phi~:~\PP^{n-1} \times \PP^{n-1} \to
\PP(M^{\mu}(n,\ \kk))$ by 
\[(a,\ b) \mapsto \left(a_i b_j + \mu_{ij} a_j b_i\right) 
\mbox {\quad for all } i,\ j = 1, \ldots , n.
\]
\end{defn}

\begin{rmk}\label{ptTomusymm}
With $a,\ b$ as in Definition~\ref{Phi}, let $q \in S_2$ be the quadratic
form 
$$q = \left(\sum_{i=1}^n a_iz_i\right)
            \left(\sum_{i=1}^n b_iz_i\right) \in \PP(S_2),$$ 
so $\mu$-rank$(q) \leq 2$. However,  using the relations of $S$, we find 
$$
q = \sum_{i=1}^n a_i b_i z_i^2 + \mathop{\sum_{i, j = 1}}_{i < j}^n 
\left(a_i b_j + \mu_{ij}a_j b_i\right)z_i z_j.
$$
It follows that $q = \tau(M)$, where $M = (a_i b_j + \mu_{ij}a_j b_i)$, 
so $M \in \kk^{\times}\Phi(a,\ b)$. Hence, $\mu$-rank$(\Phi(a,\ b))
\leq 2$ for all $a,\ b \in \PP^{n-1}$.
\end{rmk}

\begin{prop}\label{le2}
\text{\rm Im(}$\Phi$\text{\rm)} = $\{M \in \PP(M^{\mu}(n,\ \kk)) : 
\mu\mbox{\rm-rank}(M) \le 2\}$.
\end{prop}
\begin{pf}
By Remark~\ref{ptTomusymm}, Im$(\Phi) \subseteq 
\{ M \in \PP(M^{\mu}(n,\ \kk)): \mu\mbox{-rank}(M) \le 2\} = X$. 
Conversely, let $M\in X$ and write $q = \tau(M) \in \PP(S_2)$. Since 
$\mu$-rank$(q) \leq 2$, we have $q = a b$ for some $a,\ b \in 
S_1^{\times}$, where 
$a = \sum_{i=1}^n a_i z_i$ and $b = \sum_{i=1}^n b_i z_i$, with 
$a_i,\ b_i \in \kk$ for all $i$.  By Remark~\ref{ptTomusymm}, it
follows that $M = \Phi((a_i), (b_j))$.
\end{pf}

\begin{rmk}\label{xij}
Recall the notation in Remark~\ref{W}, 
and suppose $(a,\ b) \in \PP(S_1) \times \PP(S_1)$.
By our assumption in Remark~\ref{W}(a), the point 
$(a,\ b) \in \Gamma$ if and only if 
$\sum_{i, j} \alpha_{ijm} (a_i b_j + \mu_{ij}a_j b_i) = 0$ for all $m$,
where $a = (a_i)$, $b = (b_j)$; that is, if and only if the 
$\mu$-symmetric matrix $\Phi(a,\ b)$ is a zero of 
$\sum_{i, j} \alpha_{ijm} X_{ij}$ for all $m$, where $X_{ij}$ is the 
$ij$'th coordinate function on $M(n,\ \kk)$.
\end{rmk}

\begin{prop}\label{imphi}
With the assumption in Remark~\ref{W}(a),
\[\text{\rm Im(}\Phi|_{\Gamma}\text{\rm)} = 
\left\{ M \in \PP\big( \sum_{k=1}^n \kk M_k \big): 
\mu\mbox{\rm-rank}(M) \le 2 \right\}.\]
\end{prop}
\begin{pf}
Let $H = \{M \in \PP\big(\sum_{k=1}^n \kk M_k\big): 
\mu\mbox{-rank}(M) \le 2\}$ 
and let $M \in H$. Since $M$ is $\mu$-symmetric of $\mu$-rank at
most two, there exists $(a,\ b) \in \PP^{n-1} \times \PP^{n-1}$ such that
$\Phi(a,\ b) = M$, by Proposition~\ref{le2}.  Thus, by 
Lemma~\ref{qformzlocus}, $(a,\ b) \in \Gamma$, so $H \subseteq$ 
Im$(\Phi|_{\Gamma})$.

For the converse, our argument follows that of 
\cite[Proposition~1.5]{VVW}. 
Let $M = (a_{ij}) \in$ Im$(\Phi|_{\Gamma})$. So, by
Proposition~\ref{le2}, $\mu$-rank$(M) \le 2$ and, by Remark \ref{xij}, 
$\sum_{i, j=1}^n \alpha_{ijm} a_{ij} = 0$ for all $m$.
We will prove $M = \sum_{k=1}^n \beta_k M_k$, where $\beta_1,
\ldots , \beta_n \in \kk$ are defined as follows. By
Remark~\ref{yirelns}, for each $k \in \{1, \ldots , n\}$, 
$y_k \in (A_1)^2$, so $y_k = \sum_{i, j=1}^n \gamma_{ijk} Y_{ij}$, where 
$Y_{ij} = x_i x_j + \mu_{ij} x_j x_i$ and $\gamma_{ijk} \in \kk$ for all 
$i, j, k$.  For each $k = 1, \ldots , n$, we define $\beta_k \in \kk$ by 
$\beta_k = \sum_{i, j=1}^n \gamma_{ijk} a_{ij}$. 
By Remark~\ref{W}(a), 
$\sum_{i, j=1}^n \alpha_{ijm} Y_{ij} = 0$ in $A$ for all $m$, and, by 
Definition~\ref{GSCA}(a), $(Y_{ij}) = \sum_{k=1}^n M_k y_k$. Since
the behavior of the $Y_{ij}$ is mirrored by the $a_{ij}$, it follows that 
$(Y_{ij})|_{(\beta_1, \ldots , \beta_n)} = (a_{ij}) = M$, since
$ (Y_{ij})|_{(y_1, \ldots , y_n)} = (Y_{ij})$. 
Hence, 
$$ \sum_{k=1}^n \beta_k M_k = 
\sum_{k=1}^n M_k y_k |_{(\beta_1, \ldots , \beta_n )} = 
(Y_{ij})|_{(\beta_1, \ldots , \beta_n )} =M,$$
as desired. It follows that $M \in H$ and so 
Im$(\Phi|_{\Gamma}) \subseteq H$.
\end{pf}

In order to use $\Phi$ to help count the point modules over a regular 
GSCA, we need to determine which (noncommutative) quadratic forms factor 
uniquely.  To do this, we first prove, in Theorem~\ref{quadTwoWays2}, 
that a quadratic form can be factored in at most two distinct ways. 

\begin{thm}\label{quadTwoWays2}
A quadratic form can be factored in at most two distinct ways up to a
nonzero scalar multiple.
\end{thm}
\begin{pf}
Let $q \in S_2^\times$. If $q$ cannot be factored, then the result is 
trivially true. Hence, we may assume 
$$q = \left(\sum_{i=1}^n \beta_i z_i\right)\left(\sum_{i=1}^n 
\beta_i' z_i\right),$$ 
where $\beta_i,\ \beta_i' \in \kk$ for all $i$. 
If $n = 2$, then the result follows from \cite[Lemma 2.2]{VV}.
Hereafter, suppose that $n \ge 3$ and that the result holds for $n-1$ 
generators.

Case I. Suppose $\beta_i\beta_i' \ne 0$ for some $i$. Without loss of 
generality, we may assume that $i = n$ and that $\beta_n = 1 = \beta_n'$.
Suppose $q$ factors in the following three ways: 
$$
q = (a + z_n)(a' + z_n) = (b + z_n)(b' + z_n) = (c + z_n)(c' + z_n),
$$
where $a, a', b, b', c, c' \in \sum_{k=1}^{n-1} \kk z_k$.
Let $\bar{q}$ denote the image of $q$ in $S/\langle z_n \rangle$;
clearly, $\bar{q} = aa' = bb' = cc'$. The induction hypothesis implies
that $\bar{q}$ factors in at most two distinct ways up to a nonzero 
scalar multiple. Thus, 
without loss of generality, we may assume that $c = b$ and $c' = b'$. 
It follows that $q$ factors in at most two distinct ways up to a nonzero 
scalar multiple.

Case II. Suppose $\beta_i\beta_i' = 0$ for all $i$, so 
$q = \sum_{i < j} \delta_{ij} z_i z_j$ where $\delta_{ij} \in \kk$ for
all $i, j$.  We may assume, without loss of generality, that there exists
$k \in \{1, \ldots, n\}$ such that $\beta_{i} = 0$ for all $i > k$ and 
$\beta_i' = 0$ for all $i \leq k$.  By the induction hypothesis, we
may also assume that $\beta_i \neq 0$ for all $i \leq k$ and $\beta_i' 
\neq 0$ for all $i > k$.

If $q \in \langle z_i \rangle$ for some $i$, we may assume $i = n$ 
and so $k = n-1$. It follows that $q = a z_n = z_n b$, where $a, b \in 
\sum_{i=1}^{n-1} \kk z_i$. If $q = z_n b'$, where $b' \in S_1$, then 
$b = b'$ since $S$ is a domain; similarly, if $q = a' z_n$. Moreover, the
image of $q$ in the domain $S/\langle z_n \rangle$ is zero, so if also 
$q = cd$, where $c,\ d \in S_1$, then $c \in \kk z_n$ or $d \in \kk z_n$,
so $q$ factors in at most two distinct ways up to a nonzero scalar 
multiple.

Suppose $q \notin \langle z_i \rangle$ for all $i= 1,
\ldots , n$, and let $\bar{q}$ denote the image of $q$ in $S/\langle z_n 
\rangle$. By the induction hypothesis, $\bar{q}$ factors in at most
two distinct ways up to a nonzero scalar multiple, so we may assume 
$\bar{q} = ab = cd$, where $c, d \in 
\sum_{i=1}^{n-1} \kk z_i$ and $a = \sum_{i=1}^k \beta_i z_i$ and 
$b = \sum_{i=k+1}^{n-1} \beta_i' z_i$. 
Lifting to $S$, we have
$$
q = a (b + \beta_n' z_n) \quad \text{and}\quad
q = c(d+ \alpha z_n) \text{\ or\ } (c+\gamma z_n) d,
$$
where $\alpha, \gamma \in \kk^{\times}$, and these are the only ways
$q$ can factor in $S$. Hence, if $q$ factors in three distinct ways in 
$S$, then  $\beta_n' a z_n = \alpha c z_n = \gamma z_n d$, since $ab
= cd$. It follows that $c = \beta_n' \alpha^{-1} a$, since
$S$ is a domain, and $b = \beta_n' \alpha^{-1} d$, since $S/\langle z_n 
\rangle$ is a domain, and so $a (b + \beta_n' z_n)$ is a nonzero 
scalar multiple of $c(d+ \alpha z_n)$ and $\gamma$ has a unique
solution. Thus, $q$ factors in at most two distinct ways up to a 
nonzero scalar multiple.
\end{pf}

Theorem~\ref{quadTwoWays2} brings us close to our goal of generalizing 
(most of) \cite[Theorem~1.7]{VVW} from the setting of GCAs to the setting 
of GSCAs.  We first require one last technical result.

\begin{lemma} \label{rami2}
Let $\Delta_{\mu}$ denote the points $(a,\ b)\in\PP^{n-1}
\times\PP^{n-1}$ such that $(\tau \circ \Phi)(a,\ b)$ factors uniquely 
(up to nonzero scalar multiple).
The restriction of $\tau \circ \Phi$ to $(\PP^{n-1} \times \PP^{n-1}) 
\backslash \Delta_{\mu}$ has degree two and is unramified, whereas 
$\tau \circ \Phi |_{\Delta_{\mu}}$ is one-to-one.
\end{lemma}
\begin{pf}
The result is an immediate consequence of Theorem
\ref{quadTwoWays2} and the definition of $\Delta_\mu$.
\end{pf}

Our next result generalizes (most of) \cite[Theorem~1.7]{VVW}, which
we now state for comparison.

\begin{thm}\label{vvw}\cite[Theorem~1.7]{VVW}
Let~$C$ denote a GCA determined by symmetric ma\-tri\-ces $N_1, \ldots , 
N_n \in M(n,\, \kk)$ and let $\mathcal Q$ denote the corresponding 
quadric system in $\PP^{n-1}$. 
If $\mathcal Q$ has no base points, then the number 
of isomorphism classes of left (respectively, right) point modules over 
$C$ is equal to $2 r_2 + r_1 \in \N \cup \{ 0, \infty \}$, where $r_j$ 
denotes the number of matrices in $\PP\big(\sum_{k=1}^{n} \kk N_k\big)$
that have 
rank $j$.  If the number of left (respectively, right) point modules is 
finite, then $r_1 \in \{0, 1 \}$.
\end{thm}

\begin{rmk}\label{commtau}
In the setting of GCAs, if $M$ is a symmetric matrix, then $\tau(M)$
is a commutative quadratic form where $S$, in this case, is commutative;
thus, if $a,\ b \in S_1^\times$ are linearly independent, then we view
$q = ab = ba$ as two different ways to factor $q$ in $S$. It follows
that a symmetric matrix~$M$ has rank $j$, where 
$j = 1$ or $2$, if and only if $\tau(M)$ factors in $j$ distinct ways, 
up to a nonzero scalar multiple.
With this in mind, our next result is clearly a generalization of
the first part of Theorem~\ref{vvw}.
\end{rmk}

\begin{thm}\label{bigthm}
If the quadric system $\{q_1, \ldots , q_n\}$ associated to the GSCA,
$A$, is normalizing and base-point free, then the number of isomorphism 
classes of left (respectively, right) point modules over $A$ is equal to 
$2 f_2 + f_1 \in \N \cup \{0,\ \infty\}$, where $f_j$ denotes the number 
of matrices~$M$ in $\PP\big(\sum_{k=1}^n \kk M_k\big)$ such that 
$\mu$-rank$(M) \le 2$ and such that $\tau(M)$ factors in $j$ distinct 
ways (up to a nonzero scalar multiple). 
\end{thm}
\begin{pf}
Using the notation from Remark~\ref{W}, by \cite{ATV1}, the hypotheses 
on $A$ imply that the set of 
isomorphism classes of left (respectively, right) point modules over $A$ 
is in bijection with~$\Gamma$. Hence, 
the result follows from Lemma~\ref{qformzlocus}, Proposition~\ref{imphi} 
and Lemma~\ref{rami2}.
\end{pf}

The last part of Theorem~\ref{vvw} appears not to extend to the setting 
of GSCAs.  More precisely, the proof of the last part of 
Theorem~\ref{vvw} uses the correspondence between rank and factoring 
described in Remark~\ref{commtau}.  Given Remark~\ref{quadTwoWays3}, the
obvious counterpart in the setting of GSCAs is either $f_1 \in 
\{ 0,\, 1\}$ or the number of elements of $\mu$-rank one being at most 
one. However, the following two examples demonstrate that both these 
properties are unsuitable for generalizing the last part of 
Theorem~\ref{vvw} to the setting of GSCAs.

%

\begin{example}
Take $n = 4$ and let\\[-7mm]
\begin{gather*} 
\mu_{12} = \mu_{13} = \mu_{14} = -\mu_{23} = \mu_{24} = \mu_{34} = 1,
\\
q_1 = z_4^2,\qquad q_2 = z_2z_3,\qquad
q_3 = (z_1 + z_2)(z_1 + z_4),\\
q_4 = b^2 z_1^2 - a^2z_2^2 + z_3^2 + 2b z_1z_3,
\end{gather*} 
where $a,\ b \in \kk^{\times}$ and $a^2 \ne b^2$.
Since the quadric system is normalizing and base-point free, the 
corresponding GSCA, $A$, is quadratic and regular of global dimension four
(by Theorem~\ref{main2}), and is the $\kk$-algebra on generators
$x_1, \ldots, x_4$ with defining relations:\\[-4mm]
\begin{gather*} 
\begin{array}{lrcl}
x_1 x_2 + x_2 x_1 \ = \ x_1^2 - b^2x_3^2, \qquad & 
x_1 x_3 +x_3 x_1 &=& 2b x_3^2,\\[2mm] 
x_1 x_4 + x_4 x_1 \ = \  x_1^2 - b^2x_3^2,&
x_3 x_4 + x_4 x_3 &=& 0, \\[2mm] 
x_2 x_4 + x_4 x_2 \ = \ x_1^2 - b^2 x_3^2, & 
x_2^2 + a^2 x_3^2 &=& 0,
\end{array}
\end{gather*}\\[0mm] 
and has exactly eleven point modules.
In this example, $A$ is a GCA, but the algebra $S$ has been chosen 
to be noncommutative (via the choice of $\mu_{23}$).   Here, 
$\PP\big(\sum_{k=1}^4 \kk q_k\big)$ contains three elements
that factor uniquely, namely 
\[ q_1, \qquad q_4 + 2 a q_2 \qquad \text{and} \qquad q_4 - 2 a q_2. 
\]
(To see that $q_4 + 2 a q_2$ factors uniquely, we note that 
the only way it can factor is as 
$q_4 + 2 a q_2 = (b z_1 + \alpha z_2 + z_3)(b z_1 + \beta z_2 + z_3)$,
for some $\alpha,\ \beta \in \kk$, since its image factors uniquely 
in $S/\la z_2 \ra$; solving for $\alpha,\ \beta$ yields only one 
solution: $\alpha = a$, $\beta = -a$.  Similarly, for $q_4 - 2 a q_2$.)
Hence, $A$~has a finite number of point modules, yet $r_1 = 3 > 1$.
\end{example}

In the previous example, if, instead, one takes $\mu_{23} = 1$,  so
that $S$ is now commutative (as in \cite{VVW}), then the quadric system
contains only one element of rank one (up to nonzero scalar multiple), 
which agrees with Theorem~\ref{vvw}.

\begin{example}
For our second example, we consider a GSCA in \cite[\S5.3]{CV} with 
$n = 4$, where\\[-7mm]
\begin{gather*} 
q_1 = z_1 z_2, \qquad q_2 = z_3^2, \qquad
q_3 = z_1^2 - z_2z_4, \qquad q_4 = z_2^2 + z_4^2 - z_2z_3,\\
\mu_{23} = 1 = -\mu_{34}, \qquad (\mu_{14})^2 = \mu_{24} = -1,\qquad
\mu_{13} = -\mu_{14},
\end{gather*} 
so the quadric system is normalizing and base-point free. 
By Theorem \ref{main2}, the corresponding GSCA, $A$, is quadratic and 
regular of global dimension four, and is the $\kk$-algebra on generators
$x_1, \ldots, x_4$ with defining relations:\\[-2mm]
\[
\begin{array}{llr}
x_1 x_3 = \mu_{14} x_3 x_1, \qquad & x_3 x_4 = x_4 x_3, \qquad & 
x_2 x_3 + x_3 x_2= - x_4^2,\\[3mm]
x_1 x_4 = -\mu_{14} x_4 x_1, & \quad x_4^2 = x_2^2, & 
x_2 x_4 - x_4 x_2= - x_1^2,\\[3mm]
\end{array} 
\]
and has exactly five nonisomorphic point modules, two
of which correspond to $q_1 = z_1 z_2 = z_2 z_1$.
The other three point modules correspond to two quadratic forms in 
$\PP\big(\sum_{k=1}^4 \kk q_k\big)$
that have $\mu$-rank one, namely
\[
q_2 = z_3^2 \qquad \text{and} \qquad 
q_2 + 4 q_4 = (z_2 -\frac{z_3}{2} + z_4)^2  = 
(-z_2 + \frac{z_3}{2} + z_4)^2,
\]
where the latter quadratic form clearly factors in two distinct ways.
Hence, $A$ has a finite number of point modules even though two distinct 
elements of $\PP\big(\sum_{k=1}^4 \kk q_k\big)$ have $\mu$-rank one.
\end{example}


\bigskip
\bigskip

\raggedbottom

\end{document}